\documentclass[11pt,]{article}
\usepackage{t1enc}
\usepackage[latin1]{inputenc}
\usepackage[english]{babel}
\usepackage{amsfonts}
\usepackage{amssymb}
\usepackage{amsmath}
\usepackage{latexsym}
\usepackage{amsthm}
\usepackage{eufrak}
\newcounter{bff}
\theoremstyle{definition}
\newtheorem{theorem}[bff]{Theorem}

\newtheorem{proposition}[bff]{Proposition}
\newtheorem{lemma}[bff]{Lemma}
\newtheorem{definition}[bff]{Definition}
\newtheorem{example}[bff]{Example}

\def\A{\mathcal{A}}
\def\E{\mathcal{E}}
\def\D{\mathcal{D}}
\def\Oo{\mathcal{O}}
\def\L{\mathcal{L}}

\def\N{\mathbb{N}}
\def\Z{\mathbb{Z}}
\def\a{\mathfrak{A}}
\input xy
\xyoption{all}
\author{Toke Meier Carlsen}
\title{On $C^*$-algebras Associated with Sofic Shifts}
\date{September 2000}
\begin{document}
\maketitle
\begin{abstract}
We show that for a sofic shift $\Lambda$, Matsumoto's $C^*$-algebra $\Oo_\Lambda$ is isomorphic to the Cuntz-Krieger algebra of the left Krieger cover graph of $\Lambda$.
\end{abstract}

\section{Introduction}
In \cite{ck} Cuntz and Krieger defined the Cuntz-Krieger algebras. It is natural to see them as $C^*$-algebras associated with topological Markov shifts. In \cite{ass} Matsumoto associated to each subshift a $C^*$-algebra in such a way that if the subshift is a topological Markov shift, the Matsumoto algebra associated to it is the Cuntz-Krieger algebra associated to it. Furthermore in \cite{bowen} Matsumoto proved that for a sofic shift the associated Matsumoto algebra has the same $K_0$ and $K_1$ as the Cuntz-Krieger algebra for the left Krieger cover graph of the shift. It is therefore natural to ask whether the Matsumoto algebra associated to a sofic shift is isomorphic to the Cuntz-Krieger algebra of the left Krieger cover graph of the shift. In this paper we prove that indeed it is.

We will construct the isomorphism by using the universal properties of the Cuntz-Krieger algebra and the Matsumoto algebra to construct $*$-homomorphisms between them and then prove that these $*$-homomorphisms are each other's inverse.
\section{Preliminaries and notation}
Let $\a$ be a finite set endowed with the discrete topology. We will call this set the alphabet. Let $\a^\Z,\ \a^\N$ be the infinite product spaces $\prod_{i=-\infty}^\infty \a,\ \prod_{i=1}^\infty \a$ endowed with the product topology respectively. The transformation $\sigma$ on $\a^\Z,\ \a^\N$ given by $(\sigma(x))_i=x_{i+1},\ i\in \Z,\ \N$ is called the shift. Let $\Lambda$ be a shift invariant closed subset of $\a^\Z$. The topological dynamical system $(\Lambda,\sigma|_\Lambda)$ is called a subshift (cf. \cite{smale} and \cite{lm}). We denote $\sigma|_\Lambda$ by $\sigma$ for simplicity. A finite sequence $\mu=(\mu_1,\ldots ,\mu_k)$ of elements $\mu_j\in \a$ is called a finite word. We let $\a^{(\N)}$ be the set of all finite words.
We will follow the notation used in \cite{ass}, \cite{bowen}, \cite{relations} and \cite{stab}. That is we denote by $X_\Lambda$ the set of all right infinite sequences that appear in $\Lambda$, and we let for each $k\in \N$, $\Lambda^k$ be the set of all words with length $k$ appearing in some $x\in \Lambda$. We set $\Lambda_l=\bigcup_{k=0}^l \Lambda^k$ and $\Lambda^*=\bigcup_{k=0}^\infty \Lambda^k$, where $\Lambda ^0$ denotes the empty word $\emptyset$.

We will by $\Oo_\Lambda$ denote the $C^*$-algebra defined in \cite{ass} by Matsumoto. Then $\Oo_\Lambda$ is generated by partial isometries $S_i,\ i\in \a$. For $\mu\in \Lambda^*$ we define $S_\mu=S_{\mu_1}S_{\mu_2}\cdots S_{\mu_{|\mu|}}$. Following  \cite{ass}, \cite{bowen}, \cite{relations} and \cite{stab} we let $\A_\Lambda$ be the $C^*$-subalgebra of $\Oo_\Lambda$ generated by $S_\mu^*S_\mu,\ \mu\in \Lambda^*$, and $\D_\Lambda$ the $C^*$-subalgebra of $\Oo_\Lambda$ generated by $S_\mu S_\nu^*S_\nu S_\mu^*,\ \mu,\nu \in \Lambda^*$. For $\mu \in \Lambda^*$ we denote by $U_\mu$ the cylinder set for $\mu$:
\[U_\mu=\{x\in X_\Lambda \mid x_1=\mu_1,\ldots ,x_{|\mu|}=\mu_{|\mu|}\}.\]
We will need the following known facts about $\Oo_\Lambda$:
\begin{lemma}[Lemma 3.1 of \cite{ass}] \label{lille}
If $i,j\in \a$ are different, then \[S_i^*S_j=0.\]
\end{lemma}
We denote by $\mathcal{B}(X_\Lambda)$ the $C^*$-algebra of all bounded functions on $X_\Lambda$.

\begin{proposition}[Lemma 3.1 of \cite{relations}] \label{ttt}
The correspondence $\Phi$ defined by
\[\Phi(S_\mu S_\nu^*S_\nu S^*_\mu)=1_{U_\mu\cap\sigma^{-|\mu|}(\sigma^{|\nu|}(U_\nu))},\ \mu,\nu\in\Lambda^*\]
gives rise to an isomorphism from the commutative $C^*$-algebra $\D_\Lambda$
onto the $C^*$-subalgebra $C^*(1_{U_\mu\cap\sigma^{-|\mu|}(\sigma^{|\nu|}(U_\nu))},\
\mu,\nu\in\Lambda^*)$ of $\mathcal{B}(X_\Lambda)$. It's restriction to $\A_\Lambda$ yields an
isomorphism between $\A_\Lambda$ and  $C^*(1_{\sigma^{|\nu|}(U_\nu)},\
\nu\in\Lambda^*)$.
\end{proposition}

\begin{lemma} \label{diego}
For each $f\in C^*(1_{\sigma^{|\nu|}(U_\nu)},\
\nu\in\Lambda^*)$ and each $i\in \a$, 
\[\Phi(S_i\Phi^{-1}(f)S^*_i)=\sigma^\star(f)1_{U_i}\]
where $\sigma^\star$ is defined by $\sigma^\star(f)(x)=f(\sigma(x)),\ x\in
X_\Lambda$.
\end{lemma}

\begin{proof}
Let 
\[A=\{f\in C^*(1_{\sigma^{|\nu|}(U_\nu)},\
\nu\in\Lambda^*)\mid \forall i\in \a:
\Phi(S_i\Phi^{-1}(f)S^*_i)=\sigma^\star(f)1_{U_i}\}.\]
We want to show that $A= C^*(1_{\sigma^{|\nu|}(U_\nu)},\
\nu\in\Lambda^*)$.

It is easy to see that $A$ is a closed subset and that it is closed under addition and conjugation.

Let $f,g\in A$  and $i\in \a$. Then 

\begin{eqnarray*}
\Phi(S_i\Phi^{-1}(fg)S^*_i)&=&\Phi(S_i\Phi^{-1}(f)\Phi^{-1}(g)S^*_iS_iS^*_i)\\
&=&\Phi(S_i\Phi^{-1}(f)S^*_i)\Phi(S_i\Phi^{-1}(g)S^*_i)\\
&=&\sigma^\star(f)1_{U_i}\sigma^\star(g)1_{U_i}\\
&=&\sigma^\star(fg)1_{U_i}
\end{eqnarray*}
so $fg\in A$. Hence $A$ is also closed under multiplication. So it is a $C^*$-subalgebra of $ C^*(1_{\sigma^{|\nu|}(U_\nu)},\
\nu\in\Lambda^*)$.

Since 
\begin{eqnarray*}
\Phi(S_i \Phi^{-1}(1_{\sigma^{|\nu|}(U_\nu)})S^*_i)&=&\Phi(S_i S_\nu^*S_\nu
S^*_i)\\
&=&1_{U_i\cap\sigma^{-1}(\sigma^{|\nu|}(U_\nu))}\\
&=& 1_{U_i}1_{\sigma^{-1}(\sigma^{|\nu|}(U_\nu))}\\
&=&\sigma^\star(1_{\sigma^{|\nu|}(U_\nu)}) 1_{U_i},
\end{eqnarray*}
$1_{\sigma^{|\nu|}(U_\nu)}\in A$ for each $\nu\in\Lambda^*$.

So $A= C^*(1_{\sigma^{|\nu|}(U_\nu)},\
\nu\in\Lambda^*)$.
\end{proof}

\begin{theorem}[Theorem 4.9 of \cite{ass}\footnotemark] \label{sara}
\footnotetext{We remark that it is necessary to include all finite words in condition b) and c) to rule out the existens of a $*$-homomorphism from $\Oo_\Lambda$ to $\Oo_n$ (where $n$ is the number of letters in the alphabet) sending the generators to the generators.}
Let $A$ be a unital $C^*$-algebra. Suppose that there is a unital $*$-homomorphism $\psi$ from $\A_\Lambda$ to $A$ and there are partial isometries $s_i,\ i\in \a$ satisfying the following relations:
\begin{itemize}
\item[a)] $\sum_{i\in \a}s_is_i^*=1$.
\item[b)] $s_\mu^*s_\mu s_\nu=s_\nu s_{\mu\nu}^*s_{\mu\nu}$ for all $\mu,\nu\in \a^{(\N)}$.
\item[c)] $s_\mu^*s_\mu=\psi(S_\mu^*S_\mu)$ for all $\mu\in \a^{(\N)}$.
\end{itemize}
where $s_\mu=s_{\mu_1}\cdots s_{\mu_{|\mu|}},\ \mu=(\mu_1,\ldots ,\mu_{|\mu|})\in \Lambda^*$. Then $\psi$ extend to a unital $*$-homomorphism from $\Oo_\Lambda$ to $A$ such that $\psi(S_i)=s_i$ for all $i\in \a$.
\end{theorem}

\section{Sofic shifts}
As in \cite{bowen} and \cite{stab} we put for each $x\in X_\Lambda$ and each $l\in \N$
\[ \Lambda_l(x)=\{\mu \in \Lambda_l|\mu x\in X_\Lambda\}.\]
Two points $x,y \in X_\Lambda $ are said to be $l$-past equivalent if $\Lambda_l(x)=\Lambda_l(y)$. It is easy to see that this is an equivalence relation. We write this equivalence as $x\sim_l y$. Let $\E _i^l,\ i=1,2, \ldots m(l)$ be the set of all $l$-past equivalence classes of $X_\Lambda$. We denote by $\Omega_l=X_\Lambda /\sim_l$ the quotient space of the $l$-past equivalence classes of $X_\Lambda$.

Sofic shifts is a class of subshifts characterized by the following: A subshift $\Lambda$ is sofic if and only if
there exists $l\in \N$, such that $\Omega_k=\Omega_l$ for all $k\ge l$ (cf. \cite{weiss}, \cite{bowen}). In this case we will let $\Omega_\Lambda=\Omega_l$, $m_\Lambda=m(l)$ and $\E _i=\E _i^l$.

For a subset $\E \subseteq X_\Lambda$ and a $\mu\in \a^{(\N)}$ we let 
\[\mu \E =\{\mu x\in X_\Lambda|x\in \E \}.\] 
Notice that if $\mu \E _i\ne 0$ and $x\in \E _i$, then $\mu x\in X_\Lambda$.

When $\Lambda$ is a sofic shift we define the left Krieger cover graph of $\Lambda$ to be the labeled graph with vertex set $\{1,2,\ldots m_\Lambda\}$ and where there is for each vertex $i$ and each $j\in \a$ such that $j\E _i\ne \emptyset$ an edge labeled $j$ going from $k$ to $i$, where $k$ is the unique element of $\{1,2,\ldots ,m_\Lambda\}$ such that $j\E _i\subseteq \E _k$ (cf.  \cite{bowen} and \cite{krieger}). Notice that this graph is left-resolving (i.e. all edges ending at the same vertex have different labels).

For an edge $e$ we will by $s(e)$, $r(e)$ and $\L(e)$ denote the source, range and label of $e$.

We let $B_\Lambda$ be the matrix over the edge set $\mathfrak{E}_\Lambda$ defined by 
\begin{displaymath}
B_\Lambda(e,f)=\left\{ \begin{array}{ll}
1 & \textrm{if } r(e)=s(f) \\
0 & \textrm{else.}
\end{array} \right.
\end{displaymath}    
Then $B_\Lambda$ is a $\{0,1\}$-matrix with no zero-row or -column.

\begin{example} \label{malte}
Let $\a=\{0,1\}$ and $\Lambda\subseteq \a^\Z$ be the set of all sequences so that between two $1$'s there are an even number of $0$'s
Then $\Lambda$ is a sofic shift called the even shift (cf. \cite{lm}). One can show that $m_\Lambda=3$ and
\begin{eqnarray*}
\E_1&=& \{0^{2n}1x\mid n\in \N_0,\ x\in X_\Lambda\} \\
\E_2&=& \{0^{2n+1}1x\mid n\in \N_0,\ x\in X_\Lambda\} \\
\E_3&=& \{0^\infty\}.
\end{eqnarray*}
The left Krieger cover graph of $\Lambda$ is given by:

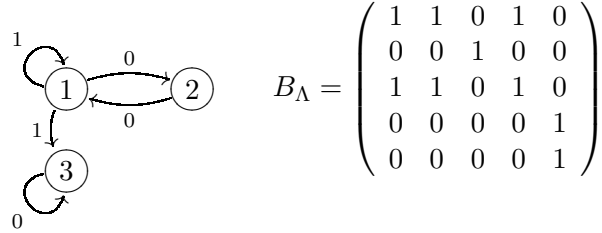
\begin{displaymath}
     \SelectTips{cm}{}
     \xymatrix @-1pc {
      *++[o][F-]{1} \ar@(l,u)[]^1 \ar@/^/[rr]^0 \ar@/_/[d]_1
      && *++[o][F-]{2} \ar@/^/[ll]^0\\
      *++[o][F-]{3} \ar@(l,d)[]_0 }
\qquad 
B_\Lambda=\left( \begin{array}{lllll}
1 & 1 & 0 & 1 & 0 \\
0 & 0 & 1 & 0 & 0 \\
1 & 1 & 0 & 1 & 0 \\
0 & 0 & 0 & 0 & 1 \\
0 & 0 & 0 & 0 & 1
\end{array} \right)
\end{displaymath}
\end{example}

The following result is a key element to the results which will be proved in this paper.

\begin{proposition} \label{vigtig}
Let $\Lambda$ be a sofic shift. Then there exists mutually orthogonal projections $E_i,\ i=1,2,\ldots ,m_\Lambda$ in $\A_\Lambda$ such that:
\begin{itemize}
\item[a)] $\A_\Lambda$ is generated by $E_i,\ i=1,2,\ldots ,m_\Lambda$.
\item[b)] For each $\mu\in \a^{(\N)}$ \[S_\mu^*S_\mu=\sum_{\mu \E _i\ne \emptyset}E_i.\]
\item[c)] For each $i\in \{1,2,\ldots m_\Lambda\}$ \[E_i=\sum_{s(e)=i}S_{\L(e)}E_{r(e)}S_{\L(e)}^*.\]
\end{itemize}
\end{proposition}

\begin{proof}
First notice that for each $\mu\in \a^{(\N)}$ 
\begin{equation} \label{sss}
\sigma^{|\mu|}(U_\mu)=\bigcup_{\mu \E _i\ne \emptyset}\E _i.
\end{equation}
From this we get that for each $i\in \{1,2,\ldots ,m_\Lambda\}$
\begin{equation*}
\E _i=\left[\bigcap_{\mu \E _i\ne \emptyset}\sigma^{|\mu|}(U_\mu)\right]\bigcap\left[\bigcap_{\mu \E _i=\emptyset}X_\Lambda\setminus \sigma^{|\mu|}(U_\mu)\right].
\end{equation*}

By (\ref{sss}) we see that there is only a finite number of different sets $\sigma^{|\mu|}(U_\mu)$. So we can for each $i\in \{1,2,\ldots m_\Lambda\}$ choose finite sets $M_i\subseteq \Lambda^*$ and $N_i\subseteq \Lambda^*$ such that

\begin{equation*}
\E _i=\left[\bigcap_{\mu\in M_i}\sigma^{|\mu|}(U_\mu)\right]\bigcap \left[\bigcap_{\mu\in N_i}X_\Lambda\setminus \sigma^{|\mu|}(U_\mu)\right].
\end{equation*}
From this we get 
\begin{equation*}
1_{\E _i}=\left[\prod_{\mu\in M_i}1_{\sigma^{|\mu|}(U_\mu)}\right] \left[\prod_{\mu\in N_i}(1-1_{\sigma^{|\mu|}(U_\mu)})\right].
\end{equation*}
So $1_{\E _i}\in C^*(1_{\sigma^{|\nu|}(U_\nu)},\ \nu\in \Lambda^*)$ for each  $i\in \{1,2,\ldots m_\Lambda\}$.

We can therefore define $E_i$ by \[E_i=\Phi^{-1}(1_{\E _i})\] where $\Phi$ is as in Proposition \ref{ttt}. Since the $\E _i$'s are mutually disjoint the $E_i$'s are mutually orthogonal projections. By (\ref{sss}) we have that
\[S_\mu^*S_\mu=\sum_{\mu \E _i\ne \emptyset}E_i.\]
So $\A_\Lambda =C^*(S_\mu^*S_\mu,\ \mu\in \Lambda^*)$ is generated by $E_i,\ i=1,2,\ldots m_\Lambda$.
Since for each $i\in \{1,2,\ldots ,m_\Lambda\}$

\begin{eqnarray*}
\E _i &=& \bigcup_{j\in \a}\{jx|jx\in \E _i\}\\
&=& \bigcup_{j\in \a}\bigcup_{k=1}^{m_\Lambda}\{jx|jx\in \E _i,x\in \E _k\}\\
&=&  \bigcup_{j\in \a}\bigcup_{j\E _k\subseteq \E _i}j\E _k\\ 
&=& \bigcup_{s(e)=i}U_{\L(e)}\cap \sigma^{-1}(\E _{r(e)})
\end{eqnarray*}
we have by Lemma \ref{diego} that 
 \[E_i=\sum_{s(e)=i}S_{\L(e)}E_{r(e)}S_{\L(e)}^*.\]
\end{proof}

\begin{example} 
If we let $\Lambda$ be as in Example \ref{malte}, then we get:
\begin{eqnarray*}
E_1 &=& S_1^*S_1(1-S_{10}^*S_{10}) \\
E_2 &=& S_{10}^*S_{10}(1-S_1^*S_1) \\
E_3 &=& S_1^*S_1S_{10}^*S_{10}.
\end{eqnarray*}
\end{example}

\section{The isomorphism}
\begin{definition}
For a matrix $A$ over a finite set $\Sigma$, with $A(i,j)\in \{0,1\}$ and where every row and column of $A$ is non-zero, we define (cf. \cite{ck}) the Cuntz-Krieger algebra for $A$ to be the universal $C^*$-algebra $\Oo_A$ generated by partial isometries $s_i,\ i\in \Sigma$ such that
\begin{itemize}
\item[a)] $s_is_i^*s_js_j^*=0$ for $i\ne j$
\item[b)] $s_i^*s_i=\sum_{j\in \Sigma}A(i,j)s_js_j^*.$
\end{itemize}
We notice that $\sum_{i\in \Sigma}s_is_i^*=1_{\Oo_A}$.
\end{definition}

\begin{proposition} \label{phi}
Let $\Lambda$ be a sofic shift. Then there exists a $*$-homomorphism from
 $\Oo_{B_\Lambda}$ to $\Oo_\Lambda$ sending $s_e$ to $S_{\L(e)}E_{r(e)}$, where $E_i$ is as in Proposition \ref{vigtig}.
\end{proposition}

\begin{proof}
Let $\widetilde{S}_e=S_{\L(e)}E_{r(e)}$. By Proposition \ref{vigtig} $S_j^*S_j=\sum_{j\E _i\ne \emptyset}E_i$ for each $j\in \a$, and since $\L(e)\E _{r(e)}\ne \emptyset$, we have $E_{r(e)}\le S_{\L(e)}^*S_{\L(e)}$. So 
\[\widetilde{S}_e^*\widetilde{S}_e=E_{r(e)} S_{\L(e)}^*S_{\L(e)}E_{r(e)}=E_{r(e)}.\]
Hence $\widetilde{S}_e$ is a partial isometry.

Since the left Krieger cover graph is left-resolving, we have that if $e\ne f$ either $\L(e)\ne \L(f)$ or $r(e)\ne r(f)$. If $\L(e)\ne \L(f)$ 
\[S_{\L(e)}^*S_{\L(f)}=0\]
by Lemma \ref{lille}, and if $r(e)\ne r(f)$ 
\[E_{r(e)}S_{\L(e)}^*S_{\L(f)}E_{r(f)}=E_{r(e)}E_{r(f)}=0.\]
So 
\[\widetilde{S}_e\widetilde{S}_e^*\widetilde{S}_f\widetilde{S}_f^*=S_{\L(e)}E_{r(e)}S_{\L(e)}^*S_{\L(f)}E_{r(f)}S_{\L(f)}^*=0\]
for $e\ne f$.

%
By Proposition \ref{vigtig}
\begin{eqnarray*}
\widetilde{S}_e^*\widetilde{S}_e &=& E_{r(e)}\\
&=& \sum_{s(f)=r(e)}S_{\L(f)}E_{r(f)}S_{\L(f)}^*\\
&=& \sum_{s(f)=r(e)}\widetilde{S}_f\widetilde{S}_f^*\\
&=& \sum_{f\in \mathfrak{E}_\Lambda} B(e,f)\widetilde{S}_f\widetilde{S}_f^*.
\end{eqnarray*}

So the partial isometries $\widetilde{S}_e,\ e\in \mathfrak{E}_\Lambda$ satisfy the Cuntz-Krieger relations and therefore there exists a $*$-homomorphism from $\Oo_{B_\Lambda}$ to $O_\Lambda$ sending $s_e$ to $\widetilde{S}_e=S_{\L(e)}E_{r(e)}$. 
\end{proof}

\begin{lemma} \label{dum}
For $\mu\in \a^{(\N)}$ and $i\in \{1,2,\ldots ,m_\Lambda\}$ the following are equivalent:
\begin{itemize}
\item[a)] $\mu \E _i\ne \emptyset$.
\item[b)] There exists a path $\alpha$ on the left Krieger cover graph of $\Lambda$ such that $\L(\alpha)=\mu$ and $r(\alpha)=i$.
\end{itemize}
The path $\alpha$ is unique, and furthermore it fulfills that $\mu \E_i\subseteq \E _{s(\alpha)}$.
\end{lemma}

\begin{proof}
We will prove the statement by induction over the length of $\mu$.
First assume that $\mu \in \a$. Then the statement follows directly from the definition of the left Krieger cover graph.

Assume next that we have proved the statement for $\mu \in \a^k$, and that $\nu \in \a^{k+1}$.
Let $\mu=(\nu_2,\nu_3,\ldots ,\nu_{|\nu|})$. 

If $\nu \E _i\ne \emptyset$, then $\mu \E _i\ne \emptyset$. So there exists a unique path $\alpha$ such that $\L(\alpha)=\mu$ and $r(\alpha)=i$ and furthermore $\mu \E _i\subseteq \E _{s(\alpha)}$. Since $\nu \E _i=\nu_1\mu\E_i\subseteq \nu_1\E _{s(\alpha)}$ and $\nu \E _i\ne \emptyset$, $\nu_1 \E _{s(\alpha)}\ne \emptyset$. Thus there exists an unique edge $e$, such that $\L(e)=\nu_1$ and $r(e)=s(\alpha)$ and furthermore $\nu_1\E _{s(\alpha)}\subseteq \E _{s(e)}$. Since $r(e)=s(\alpha)$, $e\alpha$ is a path on the left Krieger cover graph and $\L(e\alpha)=\nu$, $r(e\alpha)=i$ and $\nu \E _i\subseteq \E _{s(e\alpha)}$. 
If $\alpha'$ is another path such that $\L(\alpha')=\nu$ and $r(\alpha')=i$, then $\L((\alpha'_2,\ldots ,\alpha'_{|\alpha'|}))=\mu$, $r((\alpha'_2,\ldots ,\alpha'_{|\alpha'|}))=i$, $\L(\alpha'_1)=\nu_1$ and $r(\alpha'_1)=s((\alpha'_2,\ldots ,\alpha'_{|\alpha'|}))$. So $(\alpha'_2,\ldots ,\alpha'_{|\alpha'|})=\alpha$ and $\alpha'_1=e$. Hence $\alpha'=e\alpha$.

If there exists a path $\beta$ such that  $\L(\beta)=\nu$ and $r(\beta)=i$ 
, then $\gamma=(\beta_2,\beta_3,\ldots ,\beta_{|\beta|})$ is a path such that  $\L(\gamma)=\mu$ and $r(\beta)=i$, 
and $\beta_1$ is an edge such that  $\L(\beta_1)=\nu_1$ and $r(\beta_1)=s(\gamma)$.
So $\emptyset \ne \mu \E _i\subseteq \E _{s(\gamma)}$ and $\nu_1 \E _{s(\gamma)}\ne \emptyset$. Hence $\nu \E _i=\nu_1 \mu \E _i\ne \emptyset$.
\end{proof}
  
\begin{proposition} \label{pi}
Let $\Lambda$ be a sofic shift. Then there exists a $*$-homomorphism from
 $\Oo_{\Lambda}$ to $\Oo_{B_\Lambda}$ sending $S_i$ to $\sum_{\L(e)=i}s_e$ and $E_{r(e)}$ to $s_e^*s_e$, where $E_i$ is as in Proposition \ref{vigtig}.
\end{proposition}

\begin{proof}
Observe that $B(e,g)=B(f,g)$  for all $g\in \mathfrak{E}_\Lambda$ if $r(e)=r(f)$, and that $B(e,g)B(f,g)=0$ for all $g\in \mathfrak{E}_\Lambda$ if $r(e)\ne r(f)$. So $s_e^*s_e=s_f^*s_f$ if $r(e)=r(f)$ and $s_e^*s_es_f^*s_f=0$ if $r(e)\ne r(f)$.

Since $\A_\Lambda$ is generated by $E_i,\ i=1,2, \ldots ,m_\Lambda$ and $E_iE_j=0$ for $i\ne j$ there exists a $*$-homomorphism
$\psi$ from $\A_\Lambda$ to $\Oo_{B_\Lambda}$ sending $E_r(e)$ to $s_e^*s_e$.

For each $\mu\in \Lambda^*$ define $\widetilde{s}_\mu$ by 
\[ \widetilde{s}_\mu =\sum_{\L(\alpha)=\mu}s_{\alpha_1} \cdots s_{\alpha_{|\alpha|}}.\]
Since 
\begin{eqnarray*}
\widetilde{s}_\mu \widetilde{s}_\nu &=& \sum_{\L(\alpha)=\mu}s_{\alpha_1}\cdots s_{\alpha_{|\alpha|}}  \sum_{\L(\beta)=\nu}s_{\beta_1}\cdots s_{\beta_{|\beta|}}\\
&=& \sum_{\begin{array}{c}\scriptscriptstyle\L(\alpha)=\mu, \scriptscriptstyle\L(\beta)=\nu \\ \scriptscriptstyle r(\alpha)=s(\beta)\end{array}}s_{\alpha_1}\cdots s_{\alpha_{|\alpha|}}s_{\beta_1}\cdots s_{\beta_{|\beta|}}\\
&=&\sum_{\L(\gamma)=\mu\nu}s_{\gamma_1} \cdots s_{\gamma_{|\gamma|}}\\
&=&\widetilde{s}_{\mu\nu},
\end{eqnarray*}
we have that $\widetilde{s}_\mu=\widetilde{s}_{\mu_1}\cdots \widetilde{s}_{\mu_{|\mu|}}$ for each $\mu$.

Since the left Krieger cover graph is left-resolving we have
\begin{eqnarray*}
\widetilde{s}_i\widetilde{s}^*_i\widetilde{s}_i&=&\sum_{\L(e)=i}s_e\sum_{\L(f)=i}s_f^*\sum_{\L(g)=i}s_g\\
&=& \sum_{\L(e)=i}s_es^*_es_e\\
&=& \sum_{\L(e)=i}s_e\\
&=& \widetilde{s}_i,
\end{eqnarray*}
so $\widetilde{s}_i$ is a partial isometry.
 
We see that
\[ \sum_{i\in \a}\widetilde{s}_i\widetilde{s}^*_i=\sum_{i\in
  \a}\sum_{\L(e)=i}s_es^*_e=\sum_{e\in \mathfrak{E}_\Lambda}s_es_e^*=1.\]
If $\nu\mu\notin \Lambda^*$, 
\[\widetilde{s}^*_\nu \widetilde{s}_\nu \widetilde{s}_\mu=0=\widetilde{s}_\mu \widetilde{s}^*_{\nu\mu}\widetilde{s}_{\nu\mu},\]
and if $\nu\mu\in \Lambda^*$,  
\[\widetilde{s}^*_\nu \widetilde{s}_\nu \widetilde{s}_\mu=\widetilde{s}_\mu=\widetilde{s}_\mu \widetilde{s}^*_{\nu\mu}\widetilde{s}_{\nu\mu},\]
so  $\widetilde{s}^*_\nu \widetilde{s}_\nu \widetilde{s}_\mu=\widetilde{s}_\mu \widetilde{s}^*_{\nu\mu}\widetilde{s}_{\nu\mu}$ for all $\nu,\mu\in \a^{(\mathbb{N})}$.
%

By Proposition \ref{vigtig} and Lemma \ref{dum} we have that
\begin{eqnarray*}
\psi (S_\mu^*S_\mu)&=& \psi (\sum_{\mu \E _i\ne \emptyset}E_i)\\
&=& \sum_{\L(\alpha)=\mu}s_{\alpha_{|\alpha|}}^*s_{\alpha_{|\alpha|}}\\
&=& \sum_{\L(\alpha)=\mu}s_{\alpha_{|\alpha|}}^*\cdots s_{\alpha_1}^*s_{\alpha_1}\cdots s_{\alpha_{|\alpha|}}\\
&=& \widetilde{s}_\mu^*\widetilde{s}_\mu
\end{eqnarray*}
for all $\mu\in \a^{(\N)}$.

So according to Theorem \ref{sara} $\psi$ extends to a $*$-homomorphism from $\Oo_\Lambda$ to $\Oo_{B_\Lambda}$ sending $E_{r(e)}$ to
$s_e^*s_e$ and $S_i$ to $\widetilde{s}_i=\sum_{\L(e)=i}s_e$.
\end{proof}

\begin{theorem}
Let $\Lambda$ be a sofic shift. Then $\Oo_\Lambda \simeq \Oo_{B_\Lambda}$.
\end{theorem}

\begin{proof}
According to Proposition \ref{phi} there exists a $*$-homomorphism
 $\phi:\Oo_{B_\Lambda}\to\Oo_\Lambda$ such that
 $\phi(s_e)=S_{\L(e)}E_{r(e)}$, and according to Proposition \ref{pi} there exists a
 $*$-homomorphism $\psi:\Oo_{\Lambda}\to\Oo_{B_\Lambda}$ such that
 $\psi(S_i)=\sum_{\L(e)=i}s_e$ and $\psi(E_{r(e)})=s_e^*s_e$.

We have that 
\begin{eqnarray*}
\phi(\psi(S_i))&=&\phi(\sum_{\L(e)=i}s_e)\\
&=&\sum_{\L(e)=i}\phi(s_e)\\
&=&\sum_{\L(e)=i}S_{\L(e)}E_{r(e)}\\
&=&\sum_{\L(e)=i}S_{i}E_{r(e)}\\
&=&S_i,
\end{eqnarray*}
where we for the last equality use that $\sum_{j=1}^{m_\Lambda}E_j=1$, and that $S_iE_j=0$ if there does not exists an edge with range $j$ and label $i$. So
$\phi\circ \psi=\mathrm{id}_{\Oo_\Lambda}$, and since
\begin{eqnarray*}
\psi(\phi(s_e))&=&\psi(S_{\L(e)}E_{r(e)})\\
&=&\sum_{\L(f)=\L(e)}s_fs_e^*s_e\\
&=& s_e,
\end{eqnarray*}
$\psi\circ \phi=\mathrm{id}_{\Oo_{B_\Lambda}}$. Thus $\psi$ and $\phi$ are each other's inverse and $\Oo_\Lambda \simeq \Oo_{B_\Lambda}$.
\end{proof}

\emph{Address:}Department of Mathematics, University of Copenhagen, Universitetsparken 5, DK-2100, Copenhagen \O , Denmark.

\emph{Email:} toke@math.ku.dk. 
\end{document}